\documentclass[11pt]{article}
\usepackage{eurosym}
\usepackage{amssymb}
\usepackage{amsmath}
\usepackage{amsfonts}
\usepackage{graphicx}
\usepackage{caption}
\usepackage{subcaption}

\begin{document}

\begin{center}
\vspace{3cm}

\textbf{An Approach for Hypersurface Family with Common Geodesic Curve in
the 4D Galilean Space }$\mathbf{G}_{4}$

\bigskip

Z\"{u}hal K\"{u}\c{c}\"{u}karslan Y\"{u}zba\c{s}\i\ and Dae Won Yoon

F\i rat University, Faculty of Science, Department of Mathematics, 23119
Elazig / Turkey.

Gyeongsang National University, Department of Mathematics Education and
RINS, Jinju 660-701, South Korea

zuhal2387@yahoo.com.tr and dwyoon@gnu.ac.kr

\bigskip
\end{center}

\textbf{Abstract: }In the present study, we derive the problem of
constructing a hypersurface family from a given isogeodesic curve in the 4D
Galilean space $\mathbf{G}_{4}.$ We obtain the hypersurface as a linear
combination of the Frenet frame in $\mathbf{G}_{4}$ and examine the
necessary and sufficient conditions for the curve as a geodesic curve$.$
Finally, some examples related to our method are given for the sake of
clarity.



\section{Introduction}

Curves on a surface which locally yield the minimal distance between any two
points are of great interest. These curves are said to be geodesics which
play an important role in differential geometry. Geodesics also are curves
along which geodesic curvature vanishes. Geodesics have been studied the
subject of many studies\ in a diversity of applications, such as the
designing industry of shoes, tent manufacturing, cutting and painting path 
\cite{brond,brys,haw}.

Generally,\ the aim of mostly studies about geodesics is to set up a family
of surfaces passing a given geodesic curve and show it as a linear
combination of the marching-scale functions and the Frenet vectors. Based on
that, there have been various researches on this subject in 3-dimensional
Euclidean and non-Euclidean space \cite{Dae,kasap,kasap2,wang,al,zuhal}.

Besides, for the differential geometry of surface and hypersurface, there
exists a rising interest in 4-dimensional space \cite{kazan,duldul}. Also,
in \cite{hyper}, Bayram and Kasap gave the hypersurfaces family from a given
common geodesic curve.

In this paper, we investigate the parametric representation of hypersurface
family passing a given isogeodesic curve, i.e., both a geodesic and a
parameter curve in 4-dimensional Galilean space $\mathbf{G}_{4}$. The
remainder of our paper is given as four sections. Firstly, we mainly give
the background. Secondly, we give the parametric representations of a
hypersurface family passing a given geodesic curve and provide the necessary
and sufficient condition for that curve as a geodesic curve on the given
hypersurface.\ Subsequently, we introduce three types of the marching-scale
functions. Finally, we give some examples and figures are plotted for the
sake of clarity of our method.

\section{Preliminaries}

The Galilean space $\mathbf{G}_{3}$ is a 3-dimensional complex projective
space $P_{3}$. The absolute figure of the Galilean space comprise of $%
\{w,f,I\}$ in which $\ w$ is the ideal (absolute) plane, $\ f$ is the line
(absolute line) in $w$ and $I$ is the fixed elliptic involution of points of 
$f$.

The analyze of mechanics of plane-parallel motions reduces to the examine of
a geometry of the 3-dimensional space with $\{x,y,t\}$, is investigated by
the motion formula in \cite{yag}. It is defined that the 4D Galilean
geometry, which examines all properties invariant under motions of objects
in the space, is even complex. In an other words, it could be considered as
the properties of 4-dimensional space with coordinates that are invariant
under the general Galilean transformations in \cite{yag}.

Let $\mathcal{\ }z=\left( z_{1},z_{2},z_{3},z_{4}\right)$  and $t=$ $\left(
t_{1},t_{2},t_{3},t_{4}\right)$  be two vectors in $\mathbf{G}_{4}.$ The
Galilean scalar product in $\mathbf{G}_{4}$ is given by

\begin{equation}
\left\langle z,t\right\rangle =\left\{ 
\begin{array}{cc}
z_{1}t_{1}, & \text{if }z_{1}\neq 0\text{ or }t_{1}\neq 0 \\ 
z_{2}t_{2}+z_{3}t_{3}+z_{4}t_{4}, & \text{if \ }z_{1}=0\text{ and }t_{1}=0%
\end{array}%
\right. .  \label{1}
\end{equation}

Let $z=$ $\left( z_{1},z_{2},z_{3},z_{4}\right)$, $t=$ $\left(
t_{1},t_{2},t_{3},t_{4}\right)$  and $u=$ $\left(
u_{1},u_{2},u_{3},u_{4}\right)$  be vectors in $\mathbf{G}_{4}$. Then the
cross product in $\mathbf{G}_{4}$ is given as follows: 
\begin{equation}
z\wedge t\wedge u=\left\vert 
\begin{array}{cccc}
0 & e_{2} & e_{3} & e_{4} \\ 
z_{1} & z_{2} & z_{3} & z_{4} \\ 
t_{1} & t_{2} & t_{3} & t_{4} \\ 
u_{1} & u_{2} & u_{3} & u_{4}%
\end{array}%
\right\vert ,  \label{2}
\end{equation}%
where\ $e_{i},$ $2\leq i\leq 4,$ are the standard basis vectors.

A curve $r:I\rightarrow \mathbf{G}_{4}$ is an arbitrary curve in $\mathbf{G}%
_{4}$ is given by 
\[
r\left( t\right) =\left( f(t),g\left( t\right) ,h\left( t\right)
,l(t)\right) ,
\]%
where $f(t),g\left( t\right) ,h\left( t\right)$  and $l(t)$ are smooth
functions on \ $I\subset\mathbb{R}$.\ Let $r$ be a curve in $\mathbf{G}_{4}$%
, parametrized by the Galilean invariant arc length $s,$ is given by 
\[
r\left( s\right) =\left( s,g\left( s\right) ,h\left( s\right) ,l(s)\right).
\]%
\textbf{\ } For the curve $r$, the Frenet vectors are given in the following
forms%
\begin{eqnarray*}
t\left( s\right) &=&r^{\prime }\left( s\right) =\left( 1,g^{\prime }\left(
s\right) ,h^{\prime }\left( s\right) ,l^{\prime }(s)\right) , \\
n\left( s\right) &=&\frac{r^{\prime \prime }\left( s\right) }{\kappa \left(
s\right) }=\frac{1}{\kappa \left( s\right) }\left( 0,g^{^{\prime \prime
}}\left( s\right) ,h^{\prime \prime }\left( s\right) ,l^{^{\prime \prime
}}(s)\right) , \\
b\left( s\right) &=&\frac{1}{\tau \left( s\right) }\left( 0,(\frac{1}{\kappa
\left( s\right) }g^{^{\prime \prime }}\left( s\right) )^{\prime },(\frac{1}{%
\kappa \left( s\right) }h^{\prime \prime }\left( s\right) )^{\prime },(\frac{%
1}{\kappa \left( s\right) }l^{^{\prime \prime }}(s))^{\prime }\right) , \\
e(s) &=&\mu t\left( s\right) \wedge n\left( s\right) \wedge b\left( s\right)
,
\end{eqnarray*}%
where $\mu $ equals $\pm 1$ such that the determinant $\left\vert
t,n,b,e\right\vert =1$ and $\kappa \left( s\right) ,\tau \left( s\right)$ 
and $\sigma (s)$\ are the first, second and third curvature of $r(s)$ which
is given by, respectively,%
\begin{eqnarray}
\kappa \left( s\right) &=&\sqrt{g^{\prime \prime }\left( s\right)
^{2}+h^{\prime \prime }\left( s\right) ^{2}+l^{\prime \prime }\left(
s\right) ^{2}},  \label{3} \\
\tau \left( s\right) &=&\sqrt{\left\langle n^{\prime }\left( s\right)
,n^{\prime }\left( s\right) \right\rangle },  \nonumber \\
\sigma \left( s\right) &=&\sqrt{\left\langle b^{\prime }\left( s\right)
,e\left( s\right) \right\rangle }.  \nonumber
\end{eqnarray}

The vectors\ $t(s), n(s), b(s)$ and $e(s)$ are called the tangent, principal
normal, first binormal, and second binormal vector of $r$, respectively.

On the other hand, Frenet formulas can be given as \cite{suh}%
\begin{eqnarray}
t^{\prime }(s) &=&\kappa (s)n(s),  \label{4} \\
n^{\prime }(s) &=&\tau (s)b(s),  \nonumber \\
b^{\prime }(s) &=&-\tau (s)n(s)+\sigma (s)e(s),  \nonumber \\
e^{\prime }(s) &=&-\sigma (s)b(s).  \nonumber
\end{eqnarray}

Let $R\left( s,\varkappa ,\rho \right)$  be a hypersurface in $\mathbf{G}%
_{4}.$ The isotropic normal vector field $\eta $ of $R$ is defined as
follows \cite{dul}%
\[
\eta \left( s,\varkappa ,\rho \right) =R_{s}\wedge R_{\varkappa }\wedge
R_{\rho },
\]%
where $R_{s}=\frac{\partial R\left( s,\varkappa ,\rho \right) }{\partial s}%
, R_{\varkappa }=\frac{\partial R\left( s,\varkappa ,\rho \right) }{\partial
\varkappa }$ and $R_{\rho }=\frac{\partial R\left( s,\varkappa ,\rho \right) 
}{\partial \rho }$.

\section{Hypersurface Family with Common Geodesic Curve}

A curve $r(s)$ on a hypersurface $R\left( s,\varkappa ,\rho \right)$  in $%
\mathbf{G}_{4}$ is said to be an isoparametric curve if it is a parameter
curve, that is, there exists a pair of parameters $\varkappa _{0}$ and $\rho
_{0}$ such that $r(s)=R(s,\varkappa _{0},\rho _{0})$. Also the curve $r(s)$ on
the hypersurface $R\left( s,\varkappa ,\rho \right)$ is geodesic iff the
principal normal vector $n(s)$ of $r(s)$ is everywhere parallel to the
isotropic normal vector $\eta \left( s,\varkappa ,\rho \right)$ of the
hypersurface $R\left( s,\varkappa ,\rho \right)$. Then, a given curve $r(s)$ 
is called an isogeodesic of the hypersurface $R$ if it is both a geodesic and
an isoparametric curve on $R.$

Let $R=R\left( s,\varkappa ,\rho \right)$ be a parametric hypersurface
through the arc-length \linebreak parametrized curve $r(s)$ in $\mathbf{G}_{4}$. The
hypersurface is defined by 
\begin{eqnarray}
R\left( s,\varkappa ,\rho \right) &=&r(s)+[\alpha \left( s,\varkappa ,\rho
\right) t\left( s\right) +\beta \left( s,\varkappa ,\rho \right) n\left(
s\right)  \label{5} \\
&&+\gamma \left( s,\varkappa ,\rho \right) b\left( s\right) +\delta \left(
s,\varkappa ,\rho \right) e\left( s\right) ],  \nonumber \\
L_{1} &\leq &s\leq L_{2},\text{ }T_{1}\leq \varkappa \leq T_{2}\text{ and }%
Q_{1}\leq \rho \leq Q_{2},  \nonumber
\end{eqnarray}%
where $\alpha \left( s,\varkappa ,\rho \right)$, $\beta \left( s,\varkappa
,\rho \right)$, $\gamma \left( s,\varkappa ,\rho \right)$  and $\delta
\left( s,\varkappa ,\rho \right)$  are smooth functions. These functions are
said to be the marching-scale functions.

Our aim is to provide necessary and sufficient conditions for the given
curve $r(s)$ to be an isogeodesic curve on a hypersurface $R=R\left(
s,\varkappa ,\rho \right)$.

Firstly, let $r(s)$ be a curve on the hypersurface $R$ in $\mathbf{G}_{4}.$
If $r(s)$ is an isoparametric curve on this surface, then a parameter $%
\varkappa _{0}\in \lbrack T_{1},T_{2}]$ and $\rho _{0}\in \lbrack
Q_{1},Q_{2}]$ should be existed such that $r(s)=R\left( s,\varkappa ,\rho
\right)$, $L_{1}\leq s\leq L_{2},$ that is, 
\begin{eqnarray}
\alpha \left( s,\varkappa _{0},\rho _{0}\right) &=&\beta \left( s,\varkappa
_{0},\rho _{0}\right) =\gamma \left( s,\varkappa _{0},\rho _{0}\right)
=\delta \left( s,\varkappa _{0},\rho _{0}\right) =0,  \label{6} \\
L_{1}&\leq &s\leq L_{2},\text{\ }\varkappa _{0}\in \lbrack T_{1},T_{2}]%
\text{ and }\rho _{0}\in \lbrack Q_{1},Q_{2}].\text{ }  \nonumber
\end{eqnarray}

Secondly, $r(s)$ on the hypersurface $R\left( s,\varkappa ,\rho \right)$  is
a geodesic \ if and only if $n(s)\left\Vert \eta \left( s,\varkappa
_{0},\rho _{0}\right) \right. $.

Now, the normal vector $\eta \left( s,\varkappa _{0},\rho _{0}\right)$  can
be found by calculating the cross product of the partial derivatives and
using 
\eqref{4}
as follows:%
\begin{eqnarray}
\frac{\partial R\left( s,\varkappa ,\rho \right) }{\partial s} &=&(1+\frac{%
\partial \alpha \left( s,\varkappa ,\rho \right) }{\partial s})t\left(
s\right) +(\alpha \left( s,\varkappa ,\rho \right) \kappa \left( s\right) +%
\frac{\partial \beta \left( s,\varkappa ,\rho \right) }{\partial s}
\label{7} \\
&&-\gamma \left( s,\varkappa ,\rho \right) \tau (s))n\left( s\right) +(\beta
\left( s,\varkappa ,\rho \right) \tau (s)+\frac{\partial \gamma \left(
s,\varkappa ,\rho \right) }{\partial s}  \nonumber \\
&&-\delta \left( s,\varkappa ,\rho \right) \sigma (s))b\left( s\right)
+(\gamma \left( s,\varkappa ,\rho \right) \sigma (s)+\frac{\partial \delta
\left( s,\varkappa ,\rho \right) }{\partial s})e\left( s\right) ,  \nonumber
\end{eqnarray}

\begin{equation}
\small
{ \frac{\partial R\left( s,\varkappa ,\rho \right) }{\partial
\varkappa }=\frac{\partial \alpha \left( s,\varkappa ,\rho \right) }{%
\partial \varkappa }t\left( s\right) +\frac{\partial \beta \left(
s,\varkappa ,\rho \right) }{\partial \varkappa }n\left( s\right) +\frac{%
\partial \gamma \left( s,\varkappa ,\rho \right) }{\partial \varkappa }%
b\left( s\right) +\frac{\partial \delta \left( s,\varkappa ,\rho \right) }{%
\partial \varkappa }e\left( s\right) ,}  \label{8}
\end{equation}%
and 
\begin{equation}
\small
{ \frac{\partial R\left( s,\varkappa ,\rho \right) }{\partial \rho }=%
\frac{\partial \alpha \left( s,\varkappa ,\rho \right) }{\partial \rho }%
t\left( s\right) +\frac{\partial \beta \left( s,\varkappa ,\rho \right) }{%
\partial \rho }n\left( s\right) +\frac{\partial \gamma \left( s,\varkappa
,\rho \right) }{\partial \rho }b\left( s\right) +\frac{\partial \delta
\left( s,\varkappa ,\rho \right) }{\partial \rho }e\left( s\right) .}
\label{9}
\end{equation}%
\textbf{Remark 3.1.} Since 
\begin{eqnarray*}
\alpha \left( s,\varkappa _{0},\rho _{0}\right) &=&\beta \left( s,\varkappa
_{0},\rho _{0}\right) =\gamma \left( s,\varkappa _{0},\rho _{0}\right)
=\delta \left( s,\varkappa _{0},\rho _{0}\right) =0, \\
L_{1} &\leq &s\leq L_{2},\text{ }\varkappa _{0}\in \lbrack T_{1},T_{2}]\text{
and }\rho _{0}\in \lbrack Q_{1},Q_{2}],\text{ }
\end{eqnarray*}%
through the arc-length parametrized curve $r(s),$ by the definition of
partial differentiation, we have 
\begin{eqnarray}
\frac{\partial \alpha \left( s,\varkappa _{0},\rho _{0}\right) }{\partial s}
&=&\frac{\partial \beta \left( s,\varkappa _{0},\rho _{0}\right) }{\partial s%
}=\frac{\partial \gamma \left( s,\varkappa _{0},\rho _{0}\right) }{\partial s%
}=\frac{\partial \delta \left( s,\varkappa _{0},\rho _{0}\right) }{\partial s%
},  \label{10} \\
L_{1} &\leq &s\leq L_{2},\text{ }\varkappa _{0}\in \lbrack T_{1},T_{2}]\text{
and }\rho _{0}\in \lbrack Q_{1},Q_{2}].\text{ }  \nonumber
\end{eqnarray}

Then, from 
\eqref{5}%
, we obtain 
\begin{eqnarray}
\eta \left( s,\varkappa _{0},\rho _{0}\right) &=&\frac{\partial R\left(
s,\varkappa _{0},\rho _{0}\right) }{\partial s}\wedge \frac{\partial R\left(
s,\varkappa _{0},\rho _{0}\right) }{\partial \varkappa }\wedge \frac{%
\partial R\left( s,\varkappa _{0},\rho _{0}\right) }{\partial \rho }
\label{11} \\
&=&\varphi _{1}\left( s,\varkappa _{0},\rho _{0}\right) t\left( s\right)
-\varphi _{2}\left( s,\varkappa _{0},\rho _{0}\right) n\left( s\right) 
\nonumber \\
&&+\varphi _{3}\left( s,\varkappa _{0},\rho _{0}\right) b\left( s\right)
-\varphi _{4}\left( s,\varkappa _{0},\rho _{0}\right) e\left( s\right) . 
\nonumber
\end{eqnarray}

We need to calculate the functions $\varphi _{i}\left( s,\varkappa _{0},\rho
_{0}\right)$, $1\leq i\leq 4.$

Using 
\eqref{5}
and taking account of Remark 3.1, we have%
\begin{eqnarray}
\varphi _{1}\left( s,\varkappa _{0},\rho _{0}\right) &=&0,  \label{12} \\
\varphi _{2}\left( s,\varkappa _{0},\rho _{0}\right) &=&\frac{\partial
\gamma \left( s,\varkappa ,\rho \right) }{\partial \varkappa }\frac{\partial
\delta \left( s,\varkappa ,\rho \right) }{\partial \rho }-\frac{\partial
\gamma \left( s,\varkappa ,\rho \right) }{\partial \rho }\frac{\partial
\delta \left( s,\varkappa ,\rho \right) }{\partial \varkappa },  \nonumber \\
\varphi _{3}\left( s,\varkappa _{0},\rho _{0}\right) &=&\frac{\partial \beta
\left( s,\varkappa ,\rho \right) }{\partial \varkappa }\frac{\partial \delta
\left( s,\varkappa ,\rho \right) }{\partial \rho }-\frac{\partial \beta
\left( s,\varkappa ,\rho \right) }{\partial \rho }\frac{\partial \delta
\left( s,\varkappa ,\rho \right) }{\partial \varkappa },  \nonumber \\
\varphi _{4}\left( s,\varkappa _{0},\rho _{0}\right) &=&\frac{\partial \beta
\left( s,\varkappa ,\rho \right) }{\partial \varkappa }\frac{\partial \gamma
\left( s,\varkappa ,\rho \right) }{\partial \rho }-\frac{\partial \beta
\left( s,\varkappa ,\rho \right) }{\partial \rho }\frac{\partial \gamma
\left( s,\varkappa ,\rho \right) }{\partial \varkappa }.  \nonumber
\end{eqnarray}

So, $n(s)\left\Vert \eta \left( s,\varkappa _{0},\rho _{0}\right) \right. $
if and only if 
\[
\varphi _{2}\left( s,\varkappa _{0},\rho _{0}\right) \neq 0,\varphi
_{3}\left( s,\varkappa _{0},\rho _{0}\right) =0\text{ and }\varphi
_{4}\left( s,\varkappa _{0},\rho _{0}\right) =0.
\]

Hence, the necessary and sufficient conditions for the hypersurface $R\left(
s,\varkappa ,\rho \right)$  to have the curve $r(s)$ in $\mathbf{G}_{4}$ as
an isogeodesic curve can be given with the following theorem.

\textbf{Theorem 3.2. }Let $R\left( s,\varkappa ,\rho \right)$  be a
hypersurface having a curve $r(s)$ in $\mathbf{G}_{4}$. Then $r(s)$ is an
isogeodesic curve on the hypersurface $R$ if and only if \textbf{\ }%
\[
\alpha \left( s,\varkappa _{0},\rho _{0}\right) =\beta \left( s,\varkappa
_{0},\rho _{0}\right) =\gamma \left( s,\varkappa _{0},\rho _{0}\right)
=\delta \left( s,\varkappa _{0},\rho _{0}\right) =0,
\]%
\[
\varphi _{2}\left( s,\varkappa _{0},\rho _{0}\right) \neq 0,\varphi
_{3}\left( s,\varkappa _{0},\rho _{0}\right) =0\text{ and }\varphi
_{4}\left( s,\varkappa _{0},\rho _{0}\right) =0
\]%
satisfied, where $L_{1}\leq s\leq L_{2}$,\ $\varkappa _{0}\in \lbrack
T_{1},T_{2}]$ and $\rho _{0}\in \lbrack Q_{1},Q_{2}]$.

We call the set of hypersurfaces satisfying Theorem 3.2 an isogeodesic
hypersurface family.

\section*{MARCHING-SCALE FUNCTIONS}

For $L_{1}\leq s\leq L_{2},$ $T_{1}\leq \varkappa \leq T_{2}$ and $Q_{1}\leq
\rho \leq Q_{2},$ we will define three different above mentioned types of
the marching-scale functions.

\textbf{Type A. }Let marching-scale functions be 
\begin{eqnarray}
\alpha \left( s,\varkappa ,\rho \right) &=&\lambda \left( s\right) X\left(
\varkappa ,\rho \right) ,  \label{13} \\
\beta \left( s,\varkappa ,\rho \right) &=&\mu \left( s\right) Y\left(
\varkappa ,\rho \right) ,  \nonumber \\
\gamma \left( s,\varkappa ,\rho \right) &=&\nu \left( s\right) Z\left(
\varkappa ,\rho \right) ,  \nonumber \\
\delta \left( s,\varkappa ,\rho \right) &=&\xi \left( s\right) W\left(
\varkappa ,\rho \right) ,  \nonumber
\end{eqnarray}
where $\lambda \left( s\right), \mu \left( s\right), \nu \left( s\right), \xi \left( s\right), X\left( \varkappa ,\rho \right), Y\left( \varkappa, \rho \right), Z\left( \varkappa ,\rho \right), W\left( \varkappa, \rho
\right) \in C^{1}$ and $\lambda \left( s\right), \mu \left( s\right), \nu
\left( s\right)$ and $\xi \left( s\right)$  are not identically zero.

Hence, $r(s)$ is an isogeodesic curve on $R\left( s,\varkappa ,\rho \right)$ 
if and only if {\footnotesize 
\begin{eqnarray}
X\left( \varkappa _{0},\rho _{0}\right) &=&Y\left( \varkappa _{0},\rho
_{0}\right) =Z\left( \varkappa _{0},\rho _{0}\right) =W\left( \varkappa
_{0},\rho _{0}\right),  \label{14} \\
\nu \left( s\right) &\neq &0\text{ and }\xi \left( s\right) \neq 0\text{ and 
}\frac{\partial Z\left( s,\varkappa _{0},\rho _{0}\right) }{\partial
\varkappa }\frac{\partial W\left( s,\varkappa _{0},\rho _{0}\right) }{%
\partial \rho }-\frac{\partial Z\left( s,\varkappa _{0},\rho _{0}\right) }{%
\partial \rho }\frac{\partial W\left( s,\varkappa _{0},\rho _{0}\right) }{%
\partial \varkappa }\neq 0,  \nonumber \\
\mu \left( s\right) &=&0\text{ or }\frac{\partial Y\left( s,\varkappa
_{0},\rho _{0}\right) }{\partial \varkappa }\frac{\partial W\left(
s,\varkappa _{0},\rho _{0}\right) }{\partial \rho }-\frac{\partial Y\left(
s,\varkappa _{0},\rho _{0}\right) }{\partial \rho }\frac{\partial W\left(
s,\varkappa _{0},\rho _{0}\right) }{\partial \varkappa }=0,  \nonumber \\
\mu \left( s\right) &=&0\text{ or }\frac{\partial Y\left( s,\varkappa
_{0},\rho _{0}\right) }{\partial \varkappa }\frac{\partial Z\left(
s,\varkappa _{0},\rho _{0}\right) }{\partial \rho }-\frac{\partial Y\left(
s,\varkappa _{0},\rho _{0}\right) }{\partial \rho }\frac{\partial Z\left(
s,\varkappa _{0},\rho _{0}\right) }{\partial \varkappa }=0  \nonumber
\end{eqnarray}%
} satisfied.

To simplify 
\eqref{14}%
, we can write $\nu \left( s\right) \neq 0$ and $\xi \left( s\right) \neq 0,$
\begin{equation}
\begin{array}{c}
X\left( \varkappa _{0},\rho _{0}\right) =Y\left( \varkappa _{0},\rho
_{0}\right) =Z\left( \varkappa _{0},\rho _{0}\right) =W\left( \varkappa
_{0},\rho _{0}\right) , \\ 
\frac{\partial Z\left( s,\varkappa _{0},\rho _{0}\right) }{\partial
\varkappa }\frac{\partial W\left( s,\varkappa _{0},\rho _{0}\right) }{%
\partial \rho }-\frac{\partial Z\left( s,\varkappa _{0},\rho _{0}\right) }{%
\partial \rho }\frac{\partial W\left( s,\varkappa _{0},\rho _{0}\right) }{%
\partial \varkappa }\neq 0, \\ 
\mu \left( s\right) =0\text{ or }\frac{\partial Y\left( s,\varkappa
_{0},\rho _{0}\right) }{\partial \varkappa }=\text{ }\frac{\partial Y\left(
s,\varkappa _{0},\rho _{0}\right) }{\partial \rho }=0, \\ 
\varkappa _{0}\in \lbrack T_{1},T_{2}],\text{ }\rho _{0}\in \lbrack
Q_{1},Q_{2}].%
\end{array}
\label{15}
\end{equation}

\textbf{Type B. }Let marching-scale functions be 
\begin{eqnarray}
\alpha \left( s,\varkappa ,\rho \right) &=&\lambda \left( s,\varkappa
\right) X\left( \rho \right) ,  \label{16} \\
\beta \left( s,\varkappa ,\rho \right) &=&\mu \left( s,\varkappa \right)
Y\left( \rho \right) ,  \nonumber \\
\gamma \left( s,\varkappa ,\rho \right) &=&\nu \left( s,\varkappa \right)
Z\left( \rho \right) ,  \nonumber \\
\delta \left( s,\varkappa ,\rho \right) &=&\xi \left( s,\varkappa \right)
W\left( \rho \right) ,  \nonumber
\end{eqnarray}
where $\lambda \left( s,\varkappa \right), \mu \left( s,\varkappa \right)
,\nu \left( s,\varkappa \right) ,\xi \left( s,\varkappa \right), X\left(
\rho \right), Y\left( \rho \right), Z\left( \rho \right), W\left( \rho
\right) \in C^{1}$. Thus, $r(s)$ is an isogeodesic curve on $R\left(
s,\varkappa ,\rho \right)$  if and only if%
\begin{equation}
\begin{array}{c}
\lambda \left( s,\varkappa _{0}\right) X\left( \rho _{0}\right) =\mu \left(
s,\varkappa _{0}\right) Y\left( \rho _{0}\right) =\nu \left( s,\varkappa
_{0}\right) Z\left( \rho _{0}\right) =\xi \left( s,\varkappa _{0}\right)
W\left( \rho _{0}\right) =0, \\ 
\frac{\partial \nu \left( s,\varkappa _{0}\right) }{\partial \varkappa }\xi
\left( s,\varkappa _{0}\right) Z\left( \rho _{0}\right) \frac{dW\left( \rho
_{0}\right) }{d\rho }-\nu \left( s,\varkappa _{0}\right) \frac{\partial \xi
\left( s,\varkappa _{0}\right) }{\partial \varkappa }W\left( \rho
_{0}\right) \frac{dZ\left( \rho _{0}\right) }{d\rho }\neq 0, \\ 
Y\left( \rho _{0}\right) =\mu \left( s,\varkappa _{0}\right) =0\text{ or }%
\frac{dY\left( \rho _{0}\right) }{d\rho }=Y\left( \rho _{0}\right) =0\text{
or }\frac{dY\left( \rho _{0}\right) }{d\rho }=\frac{\partial \mu \left(
s,\varkappa _{0}\right) }{\partial \varkappa }=0\text{,} \\ 
\varkappa _{0}\in \lbrack T_{1},T_{2}],\text{ }\rho _{0}\in \lbrack
Q_{1},Q_{2}]%
\end{array}
\label{17}
\end{equation}%
satisfied.

\textbf{Type C. }Let marching-scale functions be 
\begin{eqnarray}
\alpha \left( s,\varkappa ,\rho \right) &=&\lambda \left( s,\rho \right)
X\left( \varkappa \right) ,  \label{18} \\
\beta \left( s,\varkappa ,\rho \right) &=&\mu \left( s,\rho \right) Y\left(
\varkappa \right) ,  \nonumber \\
\gamma \left( s,\varkappa ,\rho \right) &=&\nu \left( s,\rho \right) Z\left(
\varkappa \right) ,  \nonumber \\
\delta \left( s,\varkappa ,\rho \right) &=&\xi \left( s,\rho \right) W\left(
\varkappa \right) ,  \nonumber
\end{eqnarray}
where $\lambda \left( s,\rho \right), \mu \left( s,\rho \right), \nu \left(
s, \rho \right), \xi \left( s,\rho \right), X\left( \varkappa \right)
, Y\left( \varkappa \right), Z\left( \varkappa \right), W\left( \varkappa
\right) \in C^{1}$. Therefore, $r(s)$ is an isogeodesic curve on $R\left(
s,\varkappa ,\rho \right)$  if and only if%
\begin{equation}
\begin{array}{c}
\lambda \left( s,\rho _{0}\right) X\left( \varkappa _{0}\right) =\mu \left(
s,\rho _{0}\right) Y\left( \varkappa _{0}\right) =\nu \left( s,\rho
_{0}\right) Z\left( \varkappa _{0}\right) =\xi \left( s,\rho _{0}\right)
W\left( \varkappa _{0}\right) =0, \\ 
\nu \left( s,\rho _{0}\right) \frac{\partial \xi \left( s,\rho _{0}\right) }{%
\partial \rho }\frac{dZ\left( \varkappa _{0}\right) }{d\varkappa }W\left(
\varkappa _{0}\right) -\frac{\partial \nu \left( s,\rho _{0}\right) }{%
\partial \rho }\xi \left( s,\rho _{0}\right) Z\left( \varkappa _{0}\right) 
\frac{dW\left( \varkappa _{0}\right) }{d\varkappa }\neq 0, \\ 
Y\left( \varkappa _{0}\right) =\mu \left( s,\rho _{0}\right) =0\text{ or }%
\frac{dY\left( \varkappa _{0}\right) }{d\varkappa }=Y\left( \varkappa
_{0}\right) =0\text{ or }\frac{dY\left( \varkappa _{0}\right) }{d\varkappa }=%
\frac{\partial \mu \left( s,\rho _{0}\right) }{\partial \rho }=0\text{,} \\ 
\varkappa _{0}\in \lbrack T_{1},T_{2}],\text{ }\rho _{0}\in \lbrack
Q_{1},Q_{2}]%
\end{array}
\label{19}
\end{equation}%
satisfied.

\textbf{Example 3.3. }Let $r(s)$ be a curve given by parametrization 
\[
r\left( s\right) =\left( s,\cos s,\sqrt{2}\sin s,\cos s\right) .
\]

It is easy to calculate that 
\begin{eqnarray*}
t &=&\left( 1,-\sin s,\sqrt{2}\cos s,-\sin s\right) , \\
n &=&\frac{1}{\sqrt{2}}\left( 0,-\cos s,-\sqrt{2}\sin s,-\cos s\right) , \\
b &=&\frac{1}{\sqrt{2}}\left( 0,\sin s,-\sqrt{2}\cos s,\sin s\right) , \\
e &=&\frac{1}{\sqrt{2}}\left( 0,-1,0,1\right) .
\end{eqnarray*}%
Now, we obtain the hypersurface family with the isogeodesic curve $r(s)$ for
three different types of the marching-scale functions.

Marching-scale functions of Type A : Let us choose 
\begin{eqnarray*}
\lambda \left( s\right) &=&\mu \left( s\right) =\nu \left( s\right) =\xi
\left( s\right) =1, \\
X\left( \varkappa ,\rho \right) &=&\rho (\varkappa -\varkappa _{0})(\rho
-\rho _{0}), \\
Y\left( \varkappa ,\rho \right) &=&0, \\
Z\left( \varkappa ,\rho \right) &=&\rho (\varkappa -\varkappa _{0}), \\
W\left( \varkappa ,\rho \right) &=&(\rho -\rho _{0}),
\end{eqnarray*}%
where $\varkappa _{0}\in \lbrack 0,1],$ $0\leq s\leq 2\pi $ and from 
\eqref{14}
we take $\rho _{0}\neq 0.$ So, we get 
\begin{eqnarray*}
\alpha \left( s,\varkappa ,\rho \right) &=&\rho (\varkappa -\varkappa
_{0})(\rho -\rho _{0}), \\
\beta \left( s,\varkappa ,\rho \right) &=&0, \\
\gamma \left( s,\varkappa ,\rho \right) &=&\rho (\varkappa -\varkappa _{0}),
\\
\delta \left( s,\varkappa ,\rho \right) &=&(\rho -\rho _{0}),
\end{eqnarray*}%
and using 
\eqref{6}
and Frenet vectors, then we get the hypersurface which is a member of
hypersurface family as follows 
\[
\small
R\left( s,\varkappa ,\rho \right) =\left( 
\begin{array}{c}
s+\rho (\varkappa -\varkappa _{0})(\rho -\rho _{0}), \\ 
\cos s-\rho (\varkappa -\varkappa _{0})(\rho -\rho _{0})\sin s+\frac{1}{%
\sqrt{2}}\rho (\varkappa -\varkappa _{0})\sin s-\frac{1}{\sqrt{2}}(\rho
-\rho _{0}), \\ 
\sqrt{2}\sin s+\sqrt{2}\rho (\varkappa -\varkappa _{0})(\rho -\rho _{0})\cos
s-\rho (\varkappa -\varkappa _{0})\cos s, \\ 
\cos s-\rho (\varkappa -\varkappa _{0})(\rho -\rho _{0})\sin s+\frac{1}{%
\sqrt{2}}\rho (\varkappa -\varkappa _{0})\sin s+\frac{1}{\sqrt{2}}(\rho
-\rho _{0})%
\end{array}%
\right) ,
\]%
where $0\leq s\leq 2\pi$, $0\leq \varkappa _{0}\leq 1$. The position of the
curve $r(s)$ can be set on the hypersurface by changing the parameters $%
\varkappa _{0}$ and $\rho _{0}$ . Let us take $\varkappa _{0}=0$ and $\rho
_{0}=\dfrac{1}{2}$. Now $r(s)$ is again an isogeodesic on the hypersurface $%
R\left( s,\varkappa ,\rho \right)$  and the equation of the hypersurface
becomes 
\[
R\left( s,\varkappa ,\rho \right) =\left( 
\begin{array}{c}
s+\rho \varkappa (\rho -\dfrac{1}{2}), \\ 
\cos s-\rho \varkappa (\rho -\dfrac{1}{2})\sin s+\frac{1}{\sqrt{2}}\rho
\varkappa \sin s-\frac{1}{\sqrt{2}}(\rho -\dfrac{1}{2}), \\ 
\sqrt{2}\sin s+\sqrt{2}\rho \varkappa (\rho -\dfrac{1}{2})\cos s-\rho
\varkappa \cos s, \\ 
\cos s-\rho \varkappa (\rho -\dfrac{1}{2})\sin s+\frac{1}{\sqrt{2}}\rho
\varkappa \sin s+\frac{1}{\sqrt{2}}(\rho -\dfrac{1}{2})%
\end{array}%
\right) .
\]

The principle step for visualization 4D is projecting (parallel or
perspective) the geometric objects in $4$-space into the $3$-space. Thus, we
yield a three-dimensional volume. Furthermore, in practice the problem of
visualizing and approximating three-dimensional data, commonly referred to
as scalar fields. The graph of a function \ $f(\mathbf{x,y,z}):U\subset 
\mathbb{R}
^{3}\rightarrow 
\mathbb{R}
,$ $U$ is open, is a special type of parametric hypersurface the
parametrization $(\mathbf{x},\mathbf{y},\mathbf{z},\mathbf{w}=f(\mathbf{x,y,z%
}))$ in 4-space. For further information about visualization of
four-dimensional space, we refer to \cite{zho,dul,hum}. So, if we (parallel)
project the hypersurface $R\left( s,\varkappa ,\rho \right)$  into the $%
\mathbf{z=}0$ subspace and setting $\varkappa =\dfrac{1}{2},$ the surface is
given by 
\[
R_\mathbf{z}\left( s,\rho \right) =\left( 
\begin{array}{c}
s+\dfrac{1}{2}\rho (\rho -\dfrac{1}{2}), \\ 
\cos s-\dfrac{1}{2}\rho (\rho -\dfrac{1}{2})\sin s+\frac{1}{2\sqrt{2}}\rho
\sin s-\frac{1}{\sqrt{2}}(\rho -\dfrac{1}{2}), \\ 
\cos s-\dfrac{1}{2}\rho (\rho -\dfrac{1}{2})\sin s+\frac{1}{2\sqrt{2}}\rho
\sin s+\frac{1}{\sqrt{2}}(\rho -\dfrac{1}{2})%
\end{array}%
\right) ,
\]%
where $0\leq s\leq 2\pi $ and $0\leq \rho \leq 1,$ in 3-space drawn in
Figure~1-Type A.

Marching-scale functions of Type B : Let us take 
\begin{eqnarray*}
\nu \left( s,\varkappa \right) &=&(s+\varkappa ),\xi \left( s,\varkappa
\right) =s(\varkappa -\varkappa _{0}), \\
X(\rho ) &=&Y(\rho )=0, \\
Z(\rho ) &=&(\rho -\rho _{0}),W(\rho )\equiv 1,
\end{eqnarray*}%
where $\varkappa _{0}\in \lbrack 0,1],$ $\rho _{0}\in \lbrack 0,1]$ and $\pi
\leq s\leq 3\pi$. Then, we obtain 
\begin{eqnarray*}
\alpha \left( s,\varkappa ,\rho \right) &=&0, \\
\beta \left( s,\varkappa ,\rho \right) &=&0, \\
\gamma \left( s,\varkappa ,\rho \right) &=&(s+\varkappa )(\rho -\rho _{0}),
\\
\delta \left( s,\varkappa ,\rho \right) &=&s(\varkappa -\varkappa _{0}),
\end{eqnarray*}%
and using 
\eqref{6}
and Frenet vectors, the hypersurface satisfies 
\[
R\left( s,\varkappa ,\rho \right) =\left( 
\begin{array}{c}
s,\cos s+\frac{1}{\sqrt{2}}(s+\varkappa )(\rho -\rho _{0})\sin s-\frac{1}{%
\sqrt{2}}s(\varkappa -\varkappa _{0}), \\ 
\sqrt{2}\sin s-(s+\varkappa )(\rho -\rho _{0})\cos s, \\ 
\cos s+\frac{1}{\sqrt{2}}(s+\varkappa )(\rho -\rho _{0})\sin s+\frac{1}{%
\sqrt{2}}s(\varkappa -\varkappa _{0})%
\end{array}%
\right) ,
\]%
where $\pi \leq s\leq 3\pi$, $0\leq \varkappa _{0}\leq 1$ and $0\leq \rho
_{0}\leq 1.$ Then, $R\left( s,\varkappa ,\rho \right)$  is a member of the
isogeodesic hypersurface family having the curve $r(s)$ as an isogeodesic.

If $\varkappa _{0}=1$ and $\rho _{0}=0$, then the hypersurface $R$ is being 
\[
R\left( s,\varkappa ,\rho \right) =\left( 
\begin{array}{c}
s,\cos s+\frac{1}{\sqrt{2}}(s+\varkappa )\rho \sin s-\frac{1}{\sqrt{2}}%
s(\varkappa -1), \\ 
\sqrt{2}\sin s-(s+\varkappa )\rho \cos s, \\ 
\cos s+\frac{1}{\sqrt{2}}(s+\varkappa )\rho \sin s+\frac{1}{\sqrt{2}}%
s(\varkappa -1)%
\end{array}%
\right) .
\]%
Thus, if we (parallel) project the hypersurface $R\left( s,\varkappa ,\rho
\right)$  into the $\mathbf{w=}0$ subspace and fixing $\rho =\dfrac{1}{8},$
the surface is given by

\begin{equation}
R_{\mathbf{w}}\left( s,\varkappa ,\dfrac{1}{8}\right) =\left( 
\begin{array}{c}
s,\cos s+\frac{1}{8\sqrt{2}}(s+\varkappa )\sin s-\frac{1}{\sqrt{2}}%
s(\varkappa -1), \\ 
\sqrt{2}\sin s-\dfrac{1}{8}(s+\varkappa )\cos s)%
\end{array}%
\right) ,  \label{21}
\end{equation}%
where $\pi \leq s\leq 3\pi$, $0\leq \varkappa \leq 1,$ in 3-space
illustrated in Figure~1-Type B.

Marching-scale functions of Type C: Consider 
\begin{eqnarray*}
\nu \left( s,\rho \right) &=&s(\rho -\rho _{0}),\xi \left( s,\rho \right)
=(s+\rho +1), \\
X(\varkappa ) &=&Y(\varkappa )=0, \\
Z(\varkappa ) &=&\varkappa ^{2},W(\varkappa )\equiv (\varkappa -\varkappa
_{0}),
\end{eqnarray*}%
where $\rho _{0}\in \lbrack 0,1]$ and $\pi \leq s\leq 3\pi $ and from 
\eqref{14}
we take $\varkappa _{0}\neq 0.$

Then, we obtain 
\begin{eqnarray*}
\alpha \left( s,\varkappa ,\rho \right) &=&0, \\
\beta \left( s,\varkappa ,\rho \right) &=&0, \\
\gamma \left( s,\varkappa ,\rho \right) &=&s(\rho -\rho _{0})\varkappa ^{2},
\\
\gamma \left( s,\varkappa ,\rho \right) &=&(s+\rho +1)(\varkappa -\varkappa
_{0}),
\end{eqnarray*}%
and using 
\eqref{6}
and Frenet vectors, the hypersurface can be found as follows: 
\[
R\left( s,\varkappa ,\rho \right) =\left( 
\begin{array}{c}
s,\cos s+\frac{1}{\sqrt{2}}s(\rho -\rho _{0})\varkappa ^{2}\sin s-\frac{1}{%
\sqrt{2}}(s+\rho +1)(\varkappa -\varkappa _{0}), \\ 
\sqrt{2}\sin s-s(\rho -\rho _{0})\varkappa ^{2}\cos s, \\ 
\cos s+\frac{1}{\sqrt{2}}s(\rho -\rho _{0})\varkappa ^{2}\sin s+\frac{1}{%
\sqrt{2}}(s+\rho +1)(\varkappa -\varkappa _{0})%
\end{array}%
\right) .
\]

Then $R\left( s,\varkappa ,\rho \right)$  is a member of the isogeodesic
hypersurface family.

\begin{figure}[tbp]
\centering
\begin{subfigure}[b]{0.3\textwidth}
        \includegraphics[width=4cm,height=6cm]{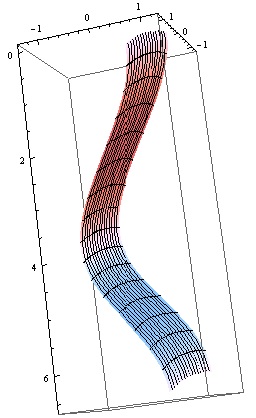}
        \caption{Type A}
        \label{fig:gull}
    \end{subfigure}
~ 
\begin{subfigure}[b]{0.3\textwidth}
        \includegraphics[width=4cm,height=6cm]{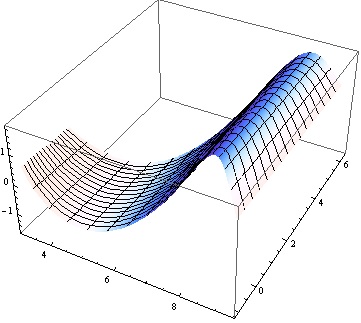}
        \caption{Type B}
        \label{fig:tiger}
    \end{subfigure}
~ 
\begin{subfigure}[b]{0.3\textwidth}
        \includegraphics[width=4cm,height=6cm]{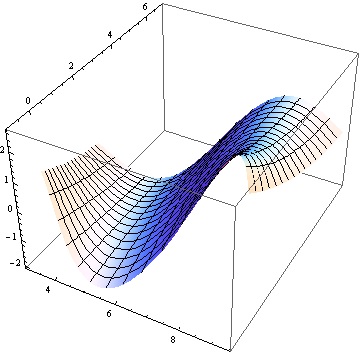}
        \caption{Type C}
        \label{fig:mouse}
    \end{subfigure}
\caption{Projection of a member of the hypersurface family with
marching-scale functions and its isogeodesic. }
\label{fig:animals}
\end{figure}

Setting $\varkappa _{0}=1$ and $\rho _{0}=0$. Then, the hypersurface $R$
becomes 
\[
R\left( s,\varkappa ,\rho \right) =\left( 
\begin{array}{c}
s,\cos s+\frac{1}{\sqrt{2}}s\rho \varkappa ^{2}\sin s-\frac{1}{\sqrt{2}}%
(s+\rho +1)(\varkappa -1), \\ 
\sqrt{2}\sin s-s\rho \varkappa ^{2}\cos s, \\ 
\cos s+\frac{1}{\sqrt{2}}s\rho \varkappa ^{2}\sin s+\frac{1}{\sqrt{2}}%
(s+\rho +1)(\varkappa -1)%
\end{array}%
\right) .
\]

Hence, if we (parallel) project the hypersurface $R\left( s,\varkappa ,\rho
\right)$  into the $\mathbf{w=}0$ subspace and fixing $\rho =\dfrac{1}{4},$
the surface is given by

\begin{equation}
R_{\mathbf{w}}\left( s,\varkappa ,\dfrac{1}{4}\right) =\left( 
\begin{array}{c}
s,\cos s+\frac{1}{4\sqrt{2}}s\varkappa ^{2}\sin s-\frac{1}{\sqrt{2}}(s+\frac{%
5}{4})(\varkappa -1), \\ 
\sqrt{2}\sin s-s\frac{1}{4}\varkappa ^{2}\cos s%
\end{array}%
\right),  \label{22}
\end{equation}%
where $\pi \leq s\leq 3\pi$, $0\leq \varkappa \leq 1,$ in 3-space plotted
in Figure~1-Type C.

\textbf{Acknowledgments:} The  second author  was supported by Basic Science Research Program through the National Research Foundation of Korea (NRF) funded by the Ministry of Education (2015R1D1A1A01060046).

\end{document}